\documentstyle[12pt,amssymb,amstex]{amsart}

\numberwithin{equation}{section}

\theoremstyle{plain}
\newtheorem{thm}{Theorem}[section]

\newtheorem{lem}[thm]{Lemma}
\newtheorem{prop}[thm]{Proposition}

\theoremstyle{definition}

\theoremstyle{remark}

\numberwithin{equation}{section}

\def\bn{{\mathbb N}}

\def\bq{{\mathbb Q}}
\def\br{{\mathbb R}}

\def\bz{{\mathbb Z}}

\def\ve{\varepsilon}

\begin{document}

\title[On equation $x^q=a$]
{On equation $x^q=a$ over $\bq_p$}

\author{Farrukh Mukhamedov}
\address{Farrukh Mukhamedov\\
Department of Computational \& Theoretical Sciences \\
Faculty of Sciences, International Islamic University Malaysia\\
P.O. Box, 141, 25710, Kuantan\\
Pahang, Malaysia} \email{{\tt far75m@@yandex.ru}}

\author{Mansoor Saburov}
\address{Mansoor Saburov\\
Department of Computational \& Theoretical Sciences \\
Faculty of Science, International Islamic University Malaysia\\
P.O. Box, 141, 25710, Kuantan\\
Pahang, Malaysia} \email{{\tt msaburov@@gmail.com}}

\begin{abstract}
In this paper, we reproduced a solvability criterion of the monomial equation
$x^q=a$ over $\bq_p$ for any natural number $q$. However, this solvability criterion was given in some algebraic number theory books by means of the multiplicative group structure of the field. Here, we present some general $p-$adic analysis method which is applicable to solve polynomial equations over the $p-$adic field. As an application of the criterion, we
describe a relationship between $q$ and $p$ in which the number $-1$
is the $q$-th power of some $p-$adic number. \vskip 0.3cm
\noindent {\it Mathematics Subject Classification 2010}: 11D88, 11S05 

\end{abstract}

\maketitle

\section{Introduction}

Over the last century, $p-$adic numbers and $p-$adic analysis have come to
play a central role in modern number theory. This importance comes from
the fact that they afford a natural and powerful language for talking about
congruences between integers, and allow the use of methods borrowed from
analysis for studying such problems.

The fields of $p-$adic numbers were introduced by German
mathematician K. Hensel \cite{Hen}. The $p-$adic numbers were
motivated primarily by an attempt to bring the ideas and techniques
of power series methods into number theory. Their canonical
representation is analogous to the expansion of analytic functions
into power series. This is one of the manifestations of the analogy
between algebraic numbers and algebraic functions.

For a fixed prime $p$, by $\bq_p$ it is denoted the field of
$p-$adic numbers, which is a completion of the rational numbers
$\bq$ with respect to the non-Archimedean norm $|\cdot|_p:\bq\to\br$
given by
\begin{eqnarray}
|x|_p=\left\{
\begin{array}{c}
  p^{-r} \ x\neq 0,\\
  0,\ \quad x=0,
\end{array}
\right.
\end{eqnarray}
here, $x=p^r\frac{m}{n}$ with $r,m\in\bz,$ $n\in\bn$,
$(m,p)=(n,p)=1$. A number $r$ is called \textit{a $p-$order} of $x$
and it is denoted by $ord_p(x)=r.$

Any $p-$adic number $x\in\bq_p$ can be uniquely represented in the
following canonical form
\begin{eqnarray*}
x=p^{ord_p(x)}\left(x_0+x_1\cdot p +x_2\cdot p^2+\cdots \right)
\end{eqnarray*}
where $x_0\in \{1,2,\cdots p-1\}$ and $x_i\in\{0,1,2,\cdots p-1\}$,
$i\geq 1,$ (see \cite{Bor Shaf}, \cite{NK}, \cite{WS})

More recently, numerous applications of $p-$adic numbers have shown
up in theoretical physics and quantum mechanics (see for example,
\cite{ArafDragFramVol}, \cite{BelGas}, \cite{FreWit}, \cite{Khren91,
Khren94}, \cite{Man}-\cite{MR2},\cite{VladVolZel,Vol}).

The $p-$adic numbers are connected with solutions of Diophantine
equations modulo increasing powers of a prime number. The study of
Diophantine equations is finding solutions of polynomial equations
or systems of equations in integers, rational numbers, or sometimes
more general number rings. Such a topic is one of the oldest
branches of number theory, in fact of mathematics itself. The theory
of Diophantine equations in number rings was well developed in
\cite{Bor Shaf}, \cite{Con}.

One of the simplest polynomial equation is the monomial equation
\begin{eqnarray}\label{eqxq=a}
x^q=a
\end{eqnarray}
where $q\in\bn$ and $a\in\bq_p$.
The solvability criterion of the monomial equation \eqref{eqxq=a} from the algebraic number theory point of view  was provided in \cite{SL}, \cite{JN}, \cite{JPS}. However, surprisingly, this solvability criterion was not mentioned in the Bible books of the $p-$adic analysis (see \cite{FGou}, \cite{NK}, \cite{WS}) except $q=2$.

Recently, J.M. Casas et.al. \cite{COR} have attempted to reproduce the solvability
criterion of the equation \eqref{eqxq=a} by concerning classification problems of high order Leibnitz algebras (see
\cite{KhudKurb}). They provided the solvability criterion in the explicit form for
two cases (i) $(q,p)=1$ and (ii) $q=p$. However,  in that paper \cite{COR}, there were some mistakes in the proof.

Let us briefly explain those mistakes.  Suppose $a,b$ are two $p-$adic numbers of the form
$$
a=a_0+a_1\cdot p+a_2\cdot p^2+\cdots,\quad \quad b=b_0+b_1\cdot p+b_2\cdot p^2+\cdots,
$$
here $a_i,b_i$ are any integer numbers with $a_0,b_0\neq 0$. In \cite{COR}, it was
implicitly stated that the numbers $a,b$ are equal each other if and only if $a_i=b_i \ (mod \ p)$ for
any $i\in\bn$. This is obviously wrong! Now we provide simple
counter examples.

Let us first consider the following two $a,b$ numbers in $\bq_3$
\begin{eqnarray*}
a=1+3+3^2+3^3+3^4+3^5+\cdots=\frac{1}{1-3}&=&-\frac12,\\
b=4+4\cdot 3^2+4\cdot 3^4+4\cdot 3^6+\cdots=4\cdot\frac{1}{1-3^2}&=&-\frac12.
\end{eqnarray*}
Then one has $a_i=1$ for any $i\geq 0$  and
$b_{2i}=4$, $b_{2i+1}=0$ for any $i\geq 0$.
However, $a_{2i+1}\neq b_{2i+1} \ (mod \ 3)$ for any $i\geq 0$. This
shows that the equality $a=b$ does not imply the relationship
$a_i=b_i \ (mod \ p)$ except $i=0$.

Now, we consider the following two $a,b$ numbers in $\bq_3$
\begin{eqnarray*}
a&=&1+3+3^2+3^3+3^4+3^5+\cdots\\
b&=&4+4\cdot 3+4\cdot 3^2+4\cdot 3^3+4\cdot 3^4+4\cdot 3^5+\cdots
\end{eqnarray*}
then we have that $a_i=1=4=b_i \ (mod \ p)$ for any $i\in\bn$.
However, $a=-\frac12$ and $b=-2$ which means $a\neq b.$ This shows
that the relationship $a_i=b_i \ (mod \ p)$ for any $i\in\bn$ does
not yield the equality $a=b$.

In this paper, we want to reproduce the solvability criterion of the monomial equation \eqref{eqxq=a} from the $p-$adic analysis point of view. Here, we shall present some general $p-$adic analysis methods which are applicable to solve some polynomial equations over $\bq_p.$ One may refer to the paper \cite{MOS} to see other applications of these methods for cubic equations over $\bq_p.$ Note that our suggested methods are completely different from the one given in \cite{COR}. As
it is usual, the solvability criterion of the
equation \eqref{eqxq=a} in $\bq_2$ and in $\bq_p$, where $p>2,$ is slightly different from each other.
Therefore, we shall separately study them. As an application, we describe the relationship between $q$
and $p$ in which the number $-1$ is the $q$-th power of some
$p-$adic number.

\section{Preliminaries}

We recall that
$$\bz_p=\{x\in\bq_{p}: |x|_p\leq1\}, \quad \bz_p^{*}=\{x\in\bq_{p}: |x|_p=1\}$$
are the set of all {\it $p-$adic integers} and {\it units} of
$\bq_p$, respectively.

Any element $x\in\bz_{p}^{*}$ has the following unique canonical form
$$x=x_0+x_1\cdot p+x_2\cdot p^2+\cdots$$
where $x_0\in \{1,2,\cdots p-1\}$ and $x_i\in\{0,1,2,\cdots p-1\}$, $i\geq 1.$

The following fact is well-known.
\begin{prop}[\cite{Bor Shaf}]
Any nonzero $p-$adic number $x$ has a unique representation of the form $x = p^{ord_p(x)}x_{*}$, where $x_{*}\in\bz_p^{*}$.
\end{prop}

We want to solve the following equation in $\bq_p$
\begin{eqnarray}\label{eq2}
x^q=a,
\end{eqnarray}
where $a\in\bq_p$. The $p-$adic numbers $x,a$ have the the following unique forms
$$x=p^{ord_p(x)}x_{*}, \quad a=p^{ord_p(a)}a_{*}.$$
After substituting these forms into \eqref{eq2} we then get that
$$p^{q\cdot ord_p(x)}x_{*}^q=p^{ord_p(a)}a_{*}$$
One can see that the equation \eqref{eq2} has a solution if
and only if
\begin{itemize}
  \item [(i)] The number $ord_p(a)$ is divided by $q$;
  \item [(ii)] The equation $x_{*}^q=a_{*}$ has a solution in $\bz_p^{*}$, whenever $a_{*}\in \bz_p^{*}.$
\end{itemize}

Therefore, it is enough to solve the equation \eqref{eq2} over
$\bz_p^{*}$ whenever $a\in \bz_p^{*}.$ Therefore, in what follows,
we shall provide a solvability criterion for the equation
\eqref{eq2} in $\bz_p^{*}$ whenever $a\in \bz_p^{*}.$

The main idea to find a solvability criterion for some polynomial equations over
$\bz_p^{*}$ is to apply Hensel's Lemma in a suitable form to the given equation.

\begin{lem}[Hensel's Lemma, \cite{Bor Shaf}]\label{Hensel}
Let $f(x)$ be polynomial whose the coefficients are $p-$adic
integers. Let $\theta$ be a $p-$adic integer such that for some
$i\geq 0$ we have
$$
f(\theta)\equiv 0 \ (mod \ p^{2i+1}),
$$
$$
f'(\theta)\equiv 0 \ (mod \ p^{i}), \quad f'(\theta)\neq 0 \ (mod \ p^{i+1}).
$$
Then the polynomial $f(x)$ has a unique $p-$adic integer root $x_0$ such that $x_0\equiv \theta\ (mod \ p^{i+1}).$
\end{lem}

\section{Auxiliary results}
In this section we are going to provide some auxiliary results which
will be used in the forthcoming sections. However, they are
independent of interest.

Throughout this paper, we call integer numbers as {\it rational
integers} in order to differentiate they from $p-$adic integers.

\begin{lem}\label{AuxiliaryFact}
Let $p$ be any prime number, $q=p^s$, $s\geq 1$, and $a\in \bz_p^{*}$. If the equation
\begin{eqnarray}\label{x^q=p^k+1}
x^{q}\equiv a\ (mod \ p^{k+1})
\end{eqnarray}
has a rational integer solution for some $k=k_0$ then  it has a
rational integer solution for any $k\geq k_0$. Here $k_0\geq s+1$
whenever $p=2$. and $k_0\geq s$ whenever $p>2$.
\end{lem}

\begin{pf}
Let $p$ be any prime number, $q=p^s$, $s\geq 1$, and $a\in
\bz_p^{*}$. We apply mathematical induction with respect to $k$
where $k\geq k_0$. Due to the assertion of the lemma, the equation
\eqref{x^q=p^k+1} has a rational integer solution for $k=k_0$.
Assume that \eqref{x^q=p^k+1} has a solution $x_n\in\bz$ for $k=n$,
where $n\geq k_0$. Now we want to show that $x_{n+1}=x_n+\ve
p^{n-s+1}$ is a solution of \eqref{x^q=p^k+1} for $k=n+1$, where
$$
\ve=a_{n+1}-\frac{x^q_n-a_0-a_1\cdot p-\cdots-a_n\cdot p^{n}}{p^{n+1}}.
$$
It is worth to mention that $\ve\in\bz,$ since $x_n$ is a solution
of \eqref{x^q=p^k+1} for $k=n$. It is clear that
\begin{eqnarray*}
x_{n+1}^q&=&x_n^q+qx_n^{q-1}\ve \cdot p^{n-s+1}+\frac{q(q-1)}{2}x_n^{q-2}\ve^2 \cdot p^{2(n-s+1)}+\cdots\\
&=&x_n^q+x_n^{q-1}\ve \cdot p^{n+1}+\frac{q-1}{2}x_n^{q-2}\ve^2\cdot p^{2n-s+2}+\cdots.
\end{eqnarray*}

If $p>2$ then $n\geq s$, $\frac{q-1}{2}\in\bn$ and $2n-s+2\geq n+2$.
If $p=2$ then $n\geq s+1$ and $2n-s+2\geq n+3$. In both cases, one
has
$$
\frac{q-1}{2}x_n^{q-2}\ve^2\cdot p^{2n-s+2}+\cdots \equiv 0 \ (mod \ p^{n+2}).
$$
Therefore, for any prime $p$, one obtain that
$$
x_{n+1}^q\equiv x_n^q+x_n^{q-1}\ve \cdot p^{n+1} \ (mod \ p^{n+2}).
$$

We know that $(x_n,p)=1$ since $a\in\bz_p^{*}$. Then one has
$x_n^{q-1}\equiv 1 \ (mod \ p)$ since $x_n^{p-1}\equiv 1 \ (mod \
p)$ and $q=p^s$. Therefore, we obtain
\begin{eqnarray*}
x_{n+1}^q&\equiv& x_n^q+x_n^{q-1}\ve\cdot p^{n+1}\equiv x_n^q+p^{n+1}\ve \ (mod \ p^{n+2})\\
&\equiv& a_0+a_1\cdot p + \cdots +a_{n+1}\cdot p^{n+1} \ (mod \ p^{n+2})\\
&\equiv&  a\ (mod \ p^{n+2}).
\end{eqnarray*}
This means that $x_{n+1}$ is a solution of \eqref{x^q=p^k+1} for
$k=n+1$. This completes the proof.
\end{pf}

\section{The solvability criterion in $\bz_p^{*}$}

In this section we are going to provide a solvability
criterion of the monomial equation \eqref{eq2} over
$\bz_p^{*},$ where $p>2$. In the sequel, we shall
suppose that $q\geq 2$, otherwise nothing to do with the equation
\eqref{eq2}.

Let $p$ be a prime number, $q\in\bn$, $a\in\bz$ with $a\neq0.$  The number $a$ is called \textit{a $q$-th power
residue modulo $p$} if the the following equation
\begin{eqnarray}\label{kthresidue}
x^q=a
\end{eqnarray} has a solution in $\bz$.

\begin{prop}[\cite{Ros}]\label{aisresidueofp}
Let $p>2$ be a prime, $q\in\bn$, $d=(q,p-1),$ and $a\in\bz$ with $a\neq0.$ Then the following statements hold true:
\begin{itemize}
  \item [(i)] $a$ is the $q$-th power residue modulo $p$ iff one has
$a^{\frac{p-1}{d}}\equiv 1 \ ( mod \ p);$
  \item [(ii)] If $a^{\frac{p-1}{d}}\equiv 1 \ ( mod \ p)$ then the equation \eqref{kthresidue} has $d$ number of solutions in $\bz$.
\end{itemize}
\end{prop}
\begin{thm}\label{Criterionforp}
Let $p>2$ be a prime and $a\in \bz_p^{*}$. The following statements hold true:
\begin{itemize}
  \item [(i)] If $(q,p)=1$ then the equation \eqref{eq2} has a solution in $\bz_p^{*}$ if and only if $a_0^{\frac{p-1}{(q,p-1)}}\equiv 1 \ ( mod \ p)$. If one has $a_0^{\frac{p-1}{(q,p-1)}}\equiv 1 \ ( mod \ p)$ then the monomial equation \eqref{eq2} has $(q,p-1)$ number of solutions in $\bz_p$.
  \item [(ii)] If $q=p^s,$ $s\geq 1$ then the equation \eqref{eq2} has a solution in $\bz_p^{*}$
  if and only if $a_0^{p^s}\equiv a \ (mod \ p^{s+1})$, i.e.,
  $$a_0^{p^s}\equiv \ a_0 + a_1\cdot p +\cdots +a_s\cdot p^s \ (mod \ p^{s+1}).$$
  Moreover, for any solution $x$ of \eqref{eq2} one has that $x\equiv a \ (mod \ p)$.
  \item [(iii)] If $q=m\cdot p^s$ with $(m,p)=1$, $s\geq 1$ then the equation \eqref{eq2}
  has a solution in $\bz_p^{*}$ if and only if $a_0^{\frac{p-1}{(m,p-1)}}\equiv 1 \ ( mod \ p)$ and $a_0^{p^s}\equiv a \ (mod \ p^{s+1})$.
\end{itemize}
\end{thm}

\begin{pf}
Let $q$ be any natural number, $p$ be an odd prime number, and $a\in
\bz_p^{*}$. Then there are three possibilities for the numbers $p$
and $q$: (i) $(q,p)=1$, (ii) $q=p^s$ with $s\geq 1$, (iii) $q=m\cdot
p^s$ with $(m,p)=1$ and $s\geq 1$.

(i). Let $(p,q)=1$.

\textsc{Only If Part.} Suppose that the equation \eqref{eq2} has a solution $x$ in $\bz_p^{*}.$ Then it is clear that
$$x^q=x^q_0+qx_0^{q-1}x_1\cdot p+\cdots=a_0+a_1\cdot p+\cdots=a.$$
This yields that $x_0^q\equiv a_0\ (mod \ p)$, i.e., $a_0$ is the $q$-th power residue modulo $p$. This means that $a_0^{\frac{p-1}{(q,p-1)}}\equiv 1 \ ( mod \ p)$.

\textsc{If Part.} Suppose that $a_0^{\frac{p-1}{(q,p-1)}}\equiv 1 \ ( mod \ p).$ Then $a_0$ is the $q$-th power residue
modulo $p$, i.e., there is $x_0$ such that $x^q_0\equiv a_0 \ (mod \
p)$.  We want to show that \eqref{eq2} has a solution in
$\bz_p^{*}$. To this end, let us consider the following function
\begin{equation}\label{f_q}
f_q(x)=x^q-a.
\end{equation}
Then, it is clear that
$$f_p(x_0)=x_0^q-a\equiv a_0-a\equiv 0 \ (mod \ p).$$ On the other hand, since $(q,p)=1$ and $(x_0,p)=1$, we get
$$f'_q(x_0)=qx_0^{q-1}\neq 0 \ (mod \ p).$$
Therefore, due to Hensel's Lemma \ref{Hensel}, we can conclude that
\eqref{eq2} has a solution in $\bz_p$. Since $x^q=a$ and $\|a\|_p=1$
for such a solution $x\in\bz_p,$ one has
$\|x\|_p=\sqrt[q]{\|a\|_q}=1$. It means that $x$ belongs to
$\bz_p^{*}.$

(ii). Let $q=p^s$ with $s\geq 1$.

\textsc{Only If Part.} Suppose that \eqref{eq2} has a solution $x$
in $\bz_p^{*}.$ Then, it is clear that
\begin{eqnarray}\label{powerofp}
x^q=x_0^q+x_0^{q-1}x_1 \cdot p^{s+1}+ p^{s+2}\left(\frac{q-1}{2}x_0^{q-2}x_1^2+\cdots\right).
\end{eqnarray}
Here, we used the fact that $\frac{q-1}{2}\in\bn$ because of $p>2.$

Therefore, one gets
\begin{eqnarray}\label{x^p=a}
x_0^q+x_0^{q-1}x_1 \cdot p^{s+1}+\cdots = a_0+a_1\cdot p+ \cdots + a_s\cdot p^s + a_{s+1}\cdot p^{s+1}+\cdots
\end{eqnarray}

This yields that $x_0^q=x_0^{p^s}\equiv a_0 \ (mod \ p).$ On the
other hand, due to Fermat's little theorem, we have $x_0^p\equiv x_0
\ (mod \ p).$ Then it follows that $a_0\equiv x_0^q\equiv x_0 \ (mod
\ p).$ We know that $x_0,a_0\in \{1,2,\cdots p-1\}$, therefore
$x_0=a_0$.

It follows again from \eqref{x^p=a} that
$$x_0^q\equiv a_0+a_1\cdot p + \cdots +a_s\cdot p^{s}\ (mod \ p^{s+1}).$$ Taking into account $x_0=a_0$ and $q=p^s$, we obtain
$$a_0^{p^s}\equiv a_0+a_1\cdot p + \cdots +a_s\cdot p^{s}\ (mod \ p^{s+1}).$$

It is worth to mention that for any solution $x\in\bz_p^{*}$ one has
$x=a \ (mod \ p)$,  because of $x_0=a_0$. In other words, the first
digit of any solution $x$ of \eqref{eq2} is the same with the first
digit of $a$.

\textsc{If Part.} Suppose that
\begin{eqnarray}\label{a_0^q=p^k+1}
a_0^{p^s}\equiv a_0+a_1\cdot p + \cdots +a_s\cdot p^{s}\ (mod \ p^{s+1}).
\end{eqnarray}
We want to show that \eqref{eq2} has a solution in $\bz_p^{*}$. In
fact, let us consider the function $f_q$ (see \eqref{f_q}).

Due to \eqref{a_0^q=p^k+1}, $a_0$ is a rational integer solution of
\eqref{x^q=p^k+1} when $k=s$. Then according to Lemma
\ref{AuxiliaryFact}, there is $x_0\in \bz$ such that
\begin{eqnarray}\label{x_0^q=p^2s+1}
x_0^{q}\equiv a\ (mod \ p^{2s+1}).
\end{eqnarray}
Then we obtain
\begin{eqnarray*}
f_q(x_0)\equiv 0 \ (mod \ p^{2s+1}).
\end{eqnarray*}

It follows from \eqref{x_0^q=p^2s+1} that $x_0^q\equiv a_0 \ (mod \
p)$ which means that $(x_0,p)=1$. We then have
$$f'_q(x_0)=qx_0^{q-1}=p^s x_0^{q-1}\equiv 0 \ (mod \ p^s),$$
$$f'_q(x_0)=qx_0^{q-1}=p^s x_0^{q-1}\neq 0 \ (mod \ p^{s+1}).$$

Therefore, due to Hensel's Lemma \ref{Hensel} we conclude that
\eqref{eq2} has a solution $x$ in $\bz_p$. From $\|a\|_p=1$, we
infer that $x$ belongs to $\bz_p^{*}.$

(iii) Let $q=m\cdot p^s$ with $(m,p)=1$ and $s\geq 1$.

\textsc{Only If Part.} Suppose that \eqref{eq2} has a solution in
$\bz_p^{*},$ i.e., there is $x\in\bz_p^{*}$ such that
$$x^q=x^{m\cdot p^s}=\left(x^{p^s}\right)^m=\left(x^m\right)^{p^s}=a.$$
This means that $a\in\bz_p^{*}$ is the $m$-th as well as the
$p^s$-th power of some $p-$adic integer numbers. Then due to the
cases (i) and (ii), we have that $a_0^{\frac{p-1}{(m,p-1)}}\equiv 1 \ ( mod \ p)$ and $a_0^{p^s}\equiv a \ (mod \ p^{s+1}),$ i.e.,
\begin{eqnarray}\label{a_0^q=p^s+1}
a_0^{p^s}\equiv \ a_0 + a_1\cdot p +\cdots +a_s\cdot p^s \ (mod \ p^{s+1}).
\end{eqnarray}

\textsc{If Part.} Suppose that $a_0$ is the $m$-th power residue
modulo $p$ with the condition \eqref{a_0^q=p^s+1}. Now we show that
\eqref{eq2} has a solution in $\bz_p^{*}$. In fact, due to the
condition \eqref{a_0^q=p^s+1} and (ii), there is $y\in\bz_p^{*}$
such that $y^{p^s}=a$ and $y\equiv a \ (mod \ p)$. This yields that
the first digit $y_0$ of the $p-$adic integer $y$ is $a_0$.
Therefore, $y_0$ is the $m$-th power residue modulo $p$. Then, due
to (i), the $p-$adic integer $y$ is $m$-th power of some $p-$adic
integer $x,$ i.e., $x^m=y.$ Consequently, the $p-$adic integer $a$
is $q$-th power of the $p-$adic integer $x,$ i.e.,
$$a=y^{p^s}=\left(x^m\right)^{p^s}=x^{m\cdot p^s}=x^q.$$

This completes the proof of the theorem.
\end{pf}

\section{The solvability criterion in $\bz_2^{*}$}

In this section we are going to provide a solvability
criterion of the monomial equation \eqref{eq2} over
$\bz_2^{*}.$  The solvability criterion in the case $p=2$ is slightly different form the case
$p>2$. Namely we have the following

\begin{thm}\label{Criterionfor2}
Let $q$ be any natural number and $a\in \bz_2^{*}$. The following statements hold true:
\begin{itemize}
  \item [(i)] If $q$ is an odd number then the equation \eqref{eq2} has a solution in $\bz_2^{*}$ for any $a$.
  \item [(ii)] If $q=2^s m $, where $m$ is an odd number and $s\geq 1$ then the equation \eqref{eq2} has a solution in $\bz_2^{*}$ if and only if $a\equiv 1 \ (mod \ 2^{s+2}),$ i.e., $$a_0=1, \quad a_1=a_2=\cdots=a_{s+1}=0.$$
  \end{itemize}
\end{thm}

\begin{pf}
Let $q$ be any natural number and $a\in \bz_2^{*}$. Then $q$ is either odd or $2^s m$, where $m$ is an odd number and $s\geq 1$.

(i) Let $q$ be an odd number. Let us show that \eqref{eq2} has a
solution in $\bz_2^{*}$ for any $a\in\bz_2^{*}$. In fact, let us
consider the same function $f_q$ as before (see \eqref{f_q}). Since
$a\in\bz_2^{*}$ and $(q,2)=1$, we have
\begin{eqnarray*}
f_q(1)&=& 1-a\equiv 0 \ (mod \ 2),\\
f'_q(1)&=& q\neq 0 \ (mod \ 2).
\end{eqnarray*}
Due to Hensel's Lemma \ref{Hensel} we can conclude that \eqref{eq2}
has a solution $x$ in $\bz_2$ and such a solution $x$ should belong
to $\bz_2^{*}$, because of $\|a\|_2=1$

(ii) We shall first consider the case $q=2^s$ with $s\geq 1$.
Suppose that \eqref{eq2} has a solution $x$ in $\bz_2^{*}.$ Then it
is clear that
\begin{eqnarray*}
x^q=x_0^q+x_0^{q-2}x_1(x_0+x_1(q-1)) \cdot 2^{s+1}+2^{s+3}\left(\frac{(q-1)(q-2)}{6}x_0^{q-2}x_1^3+\cdots\right).
\end{eqnarray*}
This yields that $x_0^q\equiv 1=a_0 \ (mod \ 2).$ Then $x_0=1$ and
\begin{eqnarray*}
x^q=1+x_1(1+x_1(q-1)) \cdot 2^{s+1}+2^{s+3}\left(\frac{(q-1)(q-2)}{6}x_0^{q-2}x_1^3+\cdots\right).
\end{eqnarray*}
Since $x_1(1+x_1(q-1))\equiv 0 (mod \ 2),$ one gets $x^q\equiv 1
(mod \ 2^{s+2}).$ Taking into account $x^{q}=a$, we heve $a\equiv 1
\ (mod \ 2^{s+2})$, i.e.,
$$a_0=1, \quad a_1=a_2=\cdots=a_{s+1}=0.$$

Let us show the converse implication. Suppose that $a\equiv 1 \ (mod
\ 2^{s+2})$. Now we show that \eqref{eq2} has a solution in
$\bz_2^{*}$ when $q=2^s$.

In fact, again consider the function $f_q$. Then $1$ is a rational
integer solution of the equation \eqref{x^q=p^k+1} with $k=s+1$ and
$p=2$, since $1\equiv a \ (mod \ 2^{s+2}).$ Then according to Lemma
\ref{AuxiliaryFact}, there is $x_0\in \bz$ such that
\begin{eqnarray*}
x_0^{q}\equiv a\ (mod \ 2^{2s+1}).
\end{eqnarray*}
We then obtain $f_q(x_0)\equiv 0 \ (mod \ p^{2s+1}).$ From
$(x_0,2)=1$, it follows  that
$$f'_q(x_0)=2^s x_0^{q-1}\equiv 0 \ (mod \ 2^s),\quad f'(x_0)=2^s x_0^{q-1}\neq 0 \ (mod \ 2^{s+1}).$$

Again Hensel's Lemma \ref{Hensel} implies that \eqref{eq2} has a
solution $x$ in $\bz_2$, which is clearly belongs to $\bz_2^{*}$.
Thus, in the case $q=2^s$, the equation \eqref{eq2} has a solution
in $\bz_2^{*}$ if and only if $a\equiv 1 \ (mod \ 2^{s+2}).$

Let us turn to the general case $q=2^s m $, where $m$ is an odd
number and $s\geq 1$. Suppose that \eqref{eq2} has a solution in
$\bz_2^{*},$ i.e., there is $x\in\bz_2^{*}$ such that
$$x^q=\left(x^m\right)^{2^s}=a.$$
This means that $a$ is $2^s$-th power of some $2-$adic integer numbers. Then, as we showed, one has $a\equiv 1 \ (mod \ 2^{s+2}).$

Now we show the reverse implication, i.e., if $a\equiv 1 \ (mod \
2^{s+2})$ then \eqref{eq2} has a solution in $\bz_2^{*}$ when $q=2^s
m$. It is clear that, since $a\equiv 1 \ (mod \ 2^{s+2})$, the
$2-$adic integer $a$ is $2^s$-th power of some $2-$adic integer
numbers $y,$ i.e., $y^{2^s}=a$. According to the case (i), since $m$
is odd, $y$ is the $m$-th power of some $2-$adic integer $x,$ i.e.,
$x^m=y.$ Hence, $a$ is the $q$-th power of the $2-$adic integer $x.$
This completes the proof of Theorem \ref{Criterionfor2}
\end{pf}

\section{When -1 is the power of some p-adic integer}

In this section, as an application of the provided criterion, we are
going to describe a relationship between $p$ and $q$ in which $-1$
is a $q$-th power of some $p-$adic integer.

\begin{thm}
Let $q$ be a natural number with $q\geq 2$ and $p$ be an odd prime number. The following statements hold true:
\begin{itemize}
  \item [(i)] The number $-1$ is any odd power of some $2-$adic integer and the number $-1$ is not any even power of any $2-$adic integer;
  \item [(ii)] If $(q,p)=1$ then the number $-1$ is a $q$-th power of some $p-$adic integer if and only if $\frac{p-1}{(q,p-1)}$ is even;
  \item [(iii)] If $q=p^s$ with $s\geq 1$ then the number $-1$ is a $q$-th power of some $p-$adic integer;
  \item [(iv)] If $q=m\cdot p^s$ with $(m,p)=1$ and $s\geq 1$ then the number $-1$ is a $q$-th power of some $p-$adic integer if and only if $\frac{p-1}{(m,p-1)}$ is even.
\end{itemize}
\end{thm}

\begin{pf}
Let $q$ be a natural number with $q\geq 2$ and $p$ be an odd prime number. It is clear that $-1$ is $2-$adic as well as $p-$adic integer. One can see that $-1$ has a canonical form with $a_i=p-1$ for any $i=0,1,2,\cdots$ in $\bq_p$. The statement (i) is obvious.

(ii) Let $(q,p)=1$. According to Theorem \ref{Criterionforp}, $-1$ is a $q$-th power of some $p-$adic integer if and only if $p-1$ is a $q$-th power residue modulo $p$. This is the same as $-1$ is a $q$-th power residue modulo $p$. Then, due to Proposition \ref{aisresidueofp}, $-1$ is a $q$-th power residue modulo $p$ if and only if $\cfrac{p-1}{(q,p-1)}$ is even.

(iii) Let $q=p^s$ with $s\geq 1$. Again according to Theorem
\ref{Criterionforp}, $-1$ is a $q$-th power of some $p-$adic integer
if and only if one has $(p-1)^q\equiv -1 \ (mod \ p^{s+1}).$ This
holds true for any $p$ and $s$. Therefore, if $q=p^s$ then $-1$ is
always $q$-th power of some $p-$adic integer.

The statement (iv) follows from Theorem \ref{Criterionforp} and
previous statements (ii), (iii). This completes the proof.
\end{pf}

\section*{Acknowledgement}  The present study have been done within
the grant FRGS0409-109 of Malaysian Ministry of Higher Education. A
part of this work was done at the Abdus Salam International Center
for Theoretical Physics (ICTP), Trieste, Italy. The author thanks
the ICTP for providing financial support during his visit as a
Junior Associate at the centre.

\end{document}